\documentclass[amstex,12pt, amssymb]{article}

\usepackage{mathtext}
\usepackage[cp1251]{inputenc}
\usepackage[T2A]{fontenc}
\usepackage[dvips]{graphicx}
\usepackage{amsmath}
\usepackage{amssymb}
\usepackage{amsxtra}
\usepackage{latexsym}
\usepackage{ifthen}

\newcounter{lemma}[section]

\newcounter{corollary}[section]

\newcounter{remark}[section]

\newcounter{theorem}[section]

\newcounter{proposition}[section]

\numberwithin{equation}{section}

\def\Xint#1{\mathchoice
   {\XXint\displaystyle\textstyle{#1}}%
   {\XXint\textstyle\scriptstyle{#1}}%
   {\XXint\scriptstyle\scriptscriptstyle{#1}}%
   {\XXint\scriptscriptstyle\scriptscriptstyle{#1}}%
   \!\int}
\def\XXint#1#2#3{{\setbox0=\hbox{$#1{#2#3}{\int}$}
     \vcenter{\hbox{$#2#3$}}\kern-.5\wd0}}
\def\dashint{\Xint-}

\def\Rn{{{\Bbb R}^n}}

\def\cc{\setcounter{equation}{0}
\setcounter{figure}{0}\setcounter{table}{0}}

\begin{document}

\markboth{\centerline{Denis Kovtonyuk, Igor Petkov  and Vladimir
Ryazanov}} {\centerline{On boundary behavior of mappings with finite
distortion in the plane}}

\author{{Denis Kovtonyuk, Igor Petkov  and Vladimir Ryazanov}}

\title{{\bf On boundary behavior of mappings \\ with finite
distortion in the plane}}

\maketitle

\large \begin{abstract} In the present paper, it was studied the
boundary behavior of the so-called lower $Q$-homeomorphisms in the
plane that are a natural generalization of the quasiconformal
mappings. In particular, it was found a series of effective
conditions on the function $Q(z)$ for a homeomorphic extension of
the given mappings to the boundary by prime ends. The
de\-ve\-lo\-ped theory is applied to mappings with finite distortion
by Iwaniec, also to solutions of the Beltrami equations, as well as
to finitely bi--Lipschitz mappings that a far--reaching extension of
the known classes of isometric and quasiisometric mappings.
\end{abstract}



\large \cc
\section{Introduction}
The problem of the boundary behavior is one of the central topics of
the theory of quasiconformal mappings and their generalizations.
During the last years, it is intensively studied various classes of
mappings with finite distortion in a natural way generalizing
conformal, quasiconformal and quasiregular mappings, see e.g. many
references in the monographs \cite{GRSY}, \cite{HK}, \cite{IM} and
\cite{MRSY}. In this case, as it was earlier, the main geometric
approach in the modern mapping theory is the method of moduli, see,
e.g., the monographs \cite{GRSY}, \cite{MRSY}, \cite{Oht},
\cite{Ri}, \cite{Vs}, \cite{Va}  and \cite{Vu}.

One of the main foundations in the present paper is the theory of
the boun\-da\-ry behavior by prime ends for the so--called lower
$Q-$homeomorphisms developed here in the first part. The second one
is Theorem 3.1 on solutions of the Beltrami equations in \cite{KPR}
that can be reformulated in the following way.

\medskip

\begin{theorem}\label{thKPR3.1} {\it Let $f$ be a homeomorphism
with finite distortion in a domain $D\subseteq{\Bbb C}$. Then $f$ is
a lower $Q$-homeomorphism at each point $z_0\in\overline{D}$ with
$Q=K_{f}$.}
\end{theorem}

\medskip

As usual, $K_{f}(z)$ denotes the {\bf dilatation} of the mapping $f$
at $z$, namely,
\begin{equation} \label{eq4.1.4} K_f(z)=
\frac{|f_z|+|f_{\overline{z}}|}{|f_z|-|f_{\overline{z}}|}
\end{equation}
if $f_z \neq 0$, $1$ if $f_z=0$ and $\infty$ otherwise, where
$f_{\overline{z}}=(f_x+if_y)/2$, $f_{z}=(f_x-if_y)/2$, $z=x+iy$, and
$f_x$ and $f_y$ are partial derivatives of $f$ in $x$ and $y$,
correspondingly.

\bigskip

Here we follow Caratheodory \cite{Car$_2$} in the definition of the
{\bf prime ends} for bounded finitely connected domains in $\Bbb C$
and refer readers to Chapter 9 in \cite{CL}, see also \cite{ABBS}
and \cite{Na} for the history of the question.

\bigskip

Reducing this case to Theorem 9.3 in \cite{CL} for simple connected
domains, we obtain the following basic fact for the theory of lower
$Q-$homeomorphisms.

\medskip

\begin{lemma}\label{thabc3} {\it Each prime end $P$ of a bounded finitely
connected domain $D$ in $\Bbb C$ contains a chain of cross--cuts
$\sigma_m$ lying on circles $S(z_0,r_m)=\{ z\in\Bbb C :\
|z-z_0|=r_m\}$ with $z_0\in\partial D$ and $r_m\to0$ as
$m\to\infty$.}
\end{lemma}

\medskip

\begin{remark}\label{METRIC}
As known, every bounded finitely connected domain $D$ in ${\Bbb C}$
can be mapped by a conformal mapping $g_0$ onto the so-called
circular domain $D_0$ whose boun\-da\-ry consists of a finite
collection of mutually disjoint circles and isolated points, see,
e.g., Theorem V.6.2 in \cite{Goluzin}. Moreover, isolated singular
points of bounded conformal mappings are removable by Theorem 1.2 in
\cite{CL} due to Weierstrass. Hence isolated points of $\partial D$
correspond to isolated points of $\partial D_0$ and inversely.


Reducing this case to the Caratheodory theorem, see, e.g., Theorem
9.4 in \cite{CL} for simple connected domains, we have a natural
one-to-one correspondence between points of $\partial D_0$ and prime
ends of the domain $D$. Denote by $\overline{D}_P$ the completion of
$D$ with its prime ends and determine in $\overline{D}_P$ the metric
$\rho_0(p_1,p_2)=\left|{\widetilde
{g_0}}(p_1)-{\widetilde{g_0}}(p_2)\right|$ where ${\widetilde
{g_0}}$ is the extension of $g_0$ to $\overline{D}_P$ mentioned
above.


If $g_*$ is another conformal mapping of the domain $D$ on a
circular domain $D_*$, then the corresponding metric
$\rho_*(p_1,p_2)=\left|{\widetilde{g_*}}(p_1)-{\widetilde{g_*}}(p_2)\right|$
generates the same convergence in $\overline{D}_P$ as the metric
$\rho_0$ because $g_0\circ g_*^{-1}$ is a conformal mapping between
the domains $D_*$ and $D_0$ that is extended to a homeomorphism
between $\overline{D_*}$ and $\overline{D_0}$. It is easy to see the
latter by applying Theorems 1.2 and 3.2 in \cite{CL}, see also Lemma
5.3 and Corollary 5.2 in \cite{IR}, correspondingly, Lemma 6.5 and
Corollary 6.12 in \cite{MRSY}. Con\-se\-quent\-ly, the given metrics
induce the same topology in the space $\overline{D}_P$ that we will
call the {\bf topology of prime ends}.


 This topology can be also described in inner terms of the domain
$D$, see, e.g., Section 9.5 in \cite{CL}, however, we prefer the
definition through the metrics because it is more clear, more
convenient and it is important for us just metrizability of
$\overline{D}_P$. Note also that the space $\overline{D}_P$ for
every bounded finitely connected domain $D$ in ${\Bbb C}$ with the
given topology is compact because the closure of the circular domain
$D_0$ is a compact space and by the construction $\widetilde
{g_0}:\overline{D}_P\to{\overline {D_0}}$ is a homeomorphism.
\end{remark}

\bigskip

Later on, we mean the continuity of mappings $f:
\overline{D}_P\to\overline{D^{\prime}}_P$ just with respect to this
topology.

\medskip

\cc
\section{On lower $Q-$homeomorphisms}

A continuous mapping $\gamma$ of an open subset $\Delta$ of the real
axis ${\Bbb R}$ or a circle into $D$ is called a {\bf dashed line},
see, e.g., Section 6.3 in \cite{MRSY}. Recall that every open set
$\Delta$ in ${\Bbb R}$ consists of a countable collection of
mutually disjoint intervals. This is the motivation for the term.

Given a family $\Gamma$ of dashed lines $\gamma$ in complex plane
${\Bbb C}$, a Borel function $\varrho:{\Bbb C}\to[0,\infty]$ is
called {\bf admissible} for $\Gamma$, write $\varrho\in{\rm
adm}\,\Gamma$, if \begin{equation}\label{eq1.2KR}
\int\limits_{\gamma}\varrho\ ds\ \geqslant\ 1\end{equation} for
every $\gamma\in\Gamma$. The {\bf (conformal) modulus} of $\Gamma$
is the quantity \begin{equation}\label{eq1.3KR}M(\Gamma)\ =\
\inf_{\varrho\in\mathrm{adm}\,\Gamma}\int\limits_{{\Bbb
C}}\varrho^2(z)\ dm(z)\end{equation} where $dm(z)$ corresponds to
the Lebesgue measure in ${\Bbb C}$. We say that a property $P$ holds
for {\bf a.e.} (almost every) $\gamma\in\Gamma$ if a subfamily of
all lines in $\Gamma$ for which $P$ fails has the modulus zero.
Later on, we also say that a Lebesgue measurable function
$\varrho:{\Bbb C}\to[0,\infty]$ is {\bf extensively admissible} for
$\Gamma$, write $\varrho\in{\rm ext\,adm}\,\Gamma$, if
(\ref{eq1.2KR}) holds for a.e. $\gamma\in\Gamma$, see, e.g., Section
9.2 in \cite{MRSY}.

The following concept was motivated by Gehring's ring definition of
qua\-si\-con\-for\-ma\-li\-ty in \cite{Ge$_1$} and first introduced
in the paper \cite{KR}. Given bounded domains $D$ and $D'$ in ${\Bbb
C}$, $z_0\in\overline{D}$, and a measurable function $Q:\Bbb
C\to(0,\infty)$, we say that a homeomorphism $f:D\to D'$ is a {\bf
lower $Q$-homeomorphism at the point} $z_0$ if
\begin{equation}\label{eq1.4KR}M(f\Sigma_{\varepsilon})\ \geqslant\ \inf\limits_{\varrho\in{\rm ext\,adm}\,\Sigma_{\varepsilon}}
\int\limits_{D\cap R_{\varepsilon}}\frac{\varrho^2(z)}{Q(z)}\
dm(z)\end{equation} for every ring
$R_{\varepsilon}=\{z\in{\overline{\Bbb
C}}:\varepsilon<|z-z_0|<\varepsilon_0\},\quad\varepsilon\in(0,\varepsilon_0),\
\varepsilon_0\in(0,d_0),$ where $d_0=\sup\limits_{z\in D}\,|z-z_0|,$
and $\Sigma_{\varepsilon}$ denotes the family of all intersections
of the circles $S(z_0,r)=\{z\in{\Bbb C}:|z-z_0|=r\},\quad
r\in(\varepsilon,\varepsilon_0),$ with the domain $D$. We also say
that a homeomorphism $f:D\to D'$ is a {\bf lower
$Q$-homeomor\-phism} if $f$ is a lower $Q$-homeomorphism at every
point $x_0\in\overline{D}$.

\medskip

Recall the criterion for homeomorphisms in ${\Bbb C}$ to be lower
$Q$-homeomorphisms, see Theorem 2.1 in \cite{KR$_1$},
correspondingly, Theorem 9.2 in \cite{MRSY}.

\medskip

\begin{proposition}\label{prOS2.2}
{\it Let $D$ and $D'$ be bounded domains in ${\Bbb C}$,
$z_0\in\overline{D}$, and $Q:\Bbb C\to(0,\infty)$ be a measurable
function. A homeomorphism $f:D\to D'$ is a lower $Q$-homeomorphism
at $z_0$ if and only if
\begin{equation}\label{eqOS2.1} M(f\Sigma_{\varepsilon})\geqslant\int\limits_{\varepsilon}^{\varepsilon_0}
\frac{dr}{||\,Q||(z_0,r)}\quad\quad\forall\
\varepsilon\in(0,\varepsilon_0)\,,\
\varepsilon_0\in(0,d_0)\,,\end{equation} where $d_0 =
\sup\limits_{z\in D}\,|\,z-z_0| $ and $||Q||(z_0,r)$ is the
$L_1$-norm of $Q$ over $D\cap S(z_0,r)$.}
\end{proposition}

\cc
\section{On regular domains}

Recall first of all the following topological notion. A domain
$D\subset{\Bbb C}$ is said to be {\bf locally connected at a point}
$z_0\in\partial D$ if, for every neighborhood $U$ of the point
$z_0$, there is a neighborhood $V\subseteq U$ of $z_0$ such that
$V\cap D$ is connected. Note that every Jordan domain $D$ in ${\Bbb
C}$ is locally connected at each point of $\partial D$, see e.g.
\cite{Wi}, p. 66.

\begin{figure}[h]
\centerline{\includegraphics[scale=0.37]{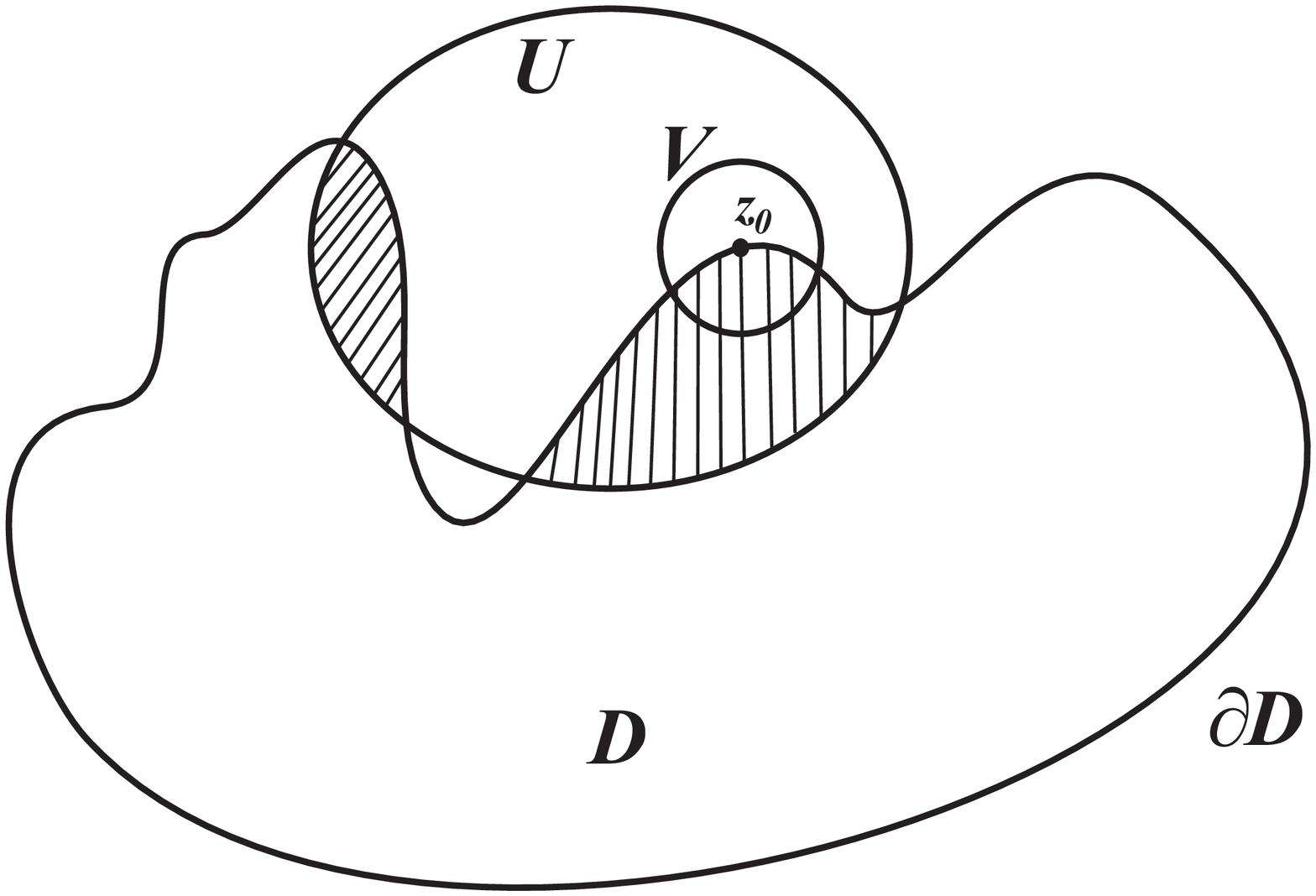}}%
\end{figure}

We say that $\partial D$ is {\bf weakly flat at a point}
$z_0\in\partial D$ if, for every neighborhood $U$ of the point $z_0$
and every number $P>0$, there is a neighborhood $V\subset U$ of
$z_0$ such that
\begin{equation}\label{eq1.5KR}M(\Delta(E,F;D))\geqslant P\end{equation} for
all continua $E$ and $F$ in $D$ intersecting $\partial U$ and
$\partial V$. Here and later on, $\Delta(E,F;D)$ denotes the family
of all paths $\gamma:[a,b]\to{\overline{\Bbb C}}$ connecting $E$ and
$F$ in $D$, i.e. $\gamma(a)\in E$, $\gamma(b)\in F$ and
$\gamma(t)\in D$ for all $t\in(a,b)$. We say that the boundary
$\partial D$ is {\bf weakly flat} if it is weakly flat at every
point in $\partial D$.

\begin{figure}[h]
\centerline{\includegraphics[scale=0.4]{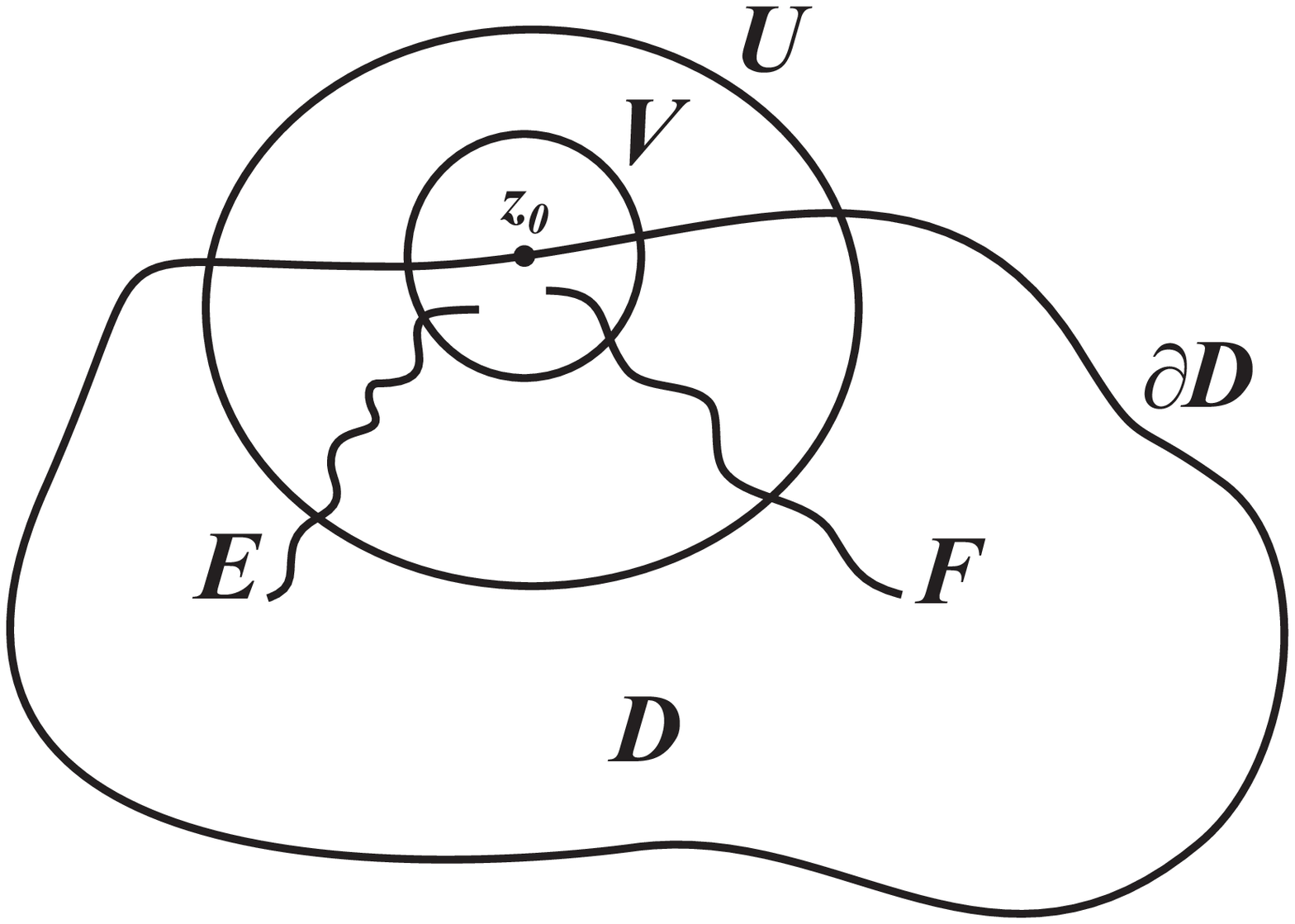}}%
\end{figure}

We also say that a point $z_0\in\partial D$ is {\bf strongly
accessible} if, for every neighborhood $U$ of the point $z_0$, there
exist a compactum $E$ in $D$, a neighborhood $V\subset U$ of $z_0$
and a number $\delta>0$ such that
\begin{equation}\label{eq1.6KR}M(\Delta(E,F;D))\geqslant\delta\end{equation} for all
continua $F$ in $D$ intersecting $\partial U$ and $\partial V$. We
say that the boundary $\partial D$ is {\bf strongly accessible} if
every point $z_0\in\partial D$ is strongly accessible.

Here, in the definitions of strongly accessible and weakly flat
boundaries, one can take as neighborhoods $U$ and $V$ of a point
$z_0$ only balls (closed or open) centered at $z_0$ or only
neighborhoods of $z_0$ in another fundamental system of
neighborhoods of $z_0$.

It is easy to see that if a domain $D$ in ${\Bbb C}$ is weakly flat
at a point $z_0\in\partial D$, then the point $z_0$ is strongly
accessible from $D$. Moreover, it was proved by us that if a domain
$D$ in ${\Bbb C}$ is weakly flat at a point $z_0\in\partial D$, then
$D$ is locally connected at $z_0$, see, e.g., Lemma 5.1 in
\cite{KR$_1$} or Lemma 3.15 in \cite{MRSY}.

The notions of strong accessibility and weak flatness at boundary
points of a domain in ${\Bbb C}$ defined in \cite{KR} are
localizations and generalizations of the cor\-res\-pon\-ding notions
introduced in \cite{MRSY$_5$}--\cite{MRSY$_6$}, cf. with the
properties $P_1$ and $P_2$ by V\"ais\"al\"a in \cite{Va} and also
with the quasiconformal accessibility and the quasiconformal
flatness by N\"akki in \cite{Na$_1$}. Many theorems on a
homeomorphic extension to the boundary of quasiconformal mappings
and their generalizations are valid under the condition of weak
flatness of boundaries. The condition of strong accessibility plays
a similar role for a continuous extension of the mappings to the
boundary. In particular, recently we have proved the following
significant statements, see either Theorem 10.1 (Lemma 6.1) in
\cite{KR$_1$} or Theorem 9.8 (Lemma 9.4) in \cite{MRSY}.

\medskip

\begin{proposition}{}\label{prKR2.1} {\it Let $D$ and $D'$ be bounded domains in ${\Bbb C}$,
$Q:D\to(0,\infty)$ a measurable function and $f:D\to D'$ a lower
$Q$-homeomorphism in $\partial D$. Suppose that the domain $D$ is
locally connected on $\partial D$ and that the domain $D'$ has a
(strongly accessible) weakly flat boundary. If
\begin{equation}\label{eqKPR2.1}\int\limits_{0}^{\delta(z_0)} \frac{dr}{||\,Q||(z_0,r)}\
=\ \infty\qquad\forall\ z_0\in\partial D\end{equation} for some
$\delta(z_0)\in(0,d(z_0))$ where $d(z_0)=\sup\limits_{z\in
D}\,|\,z-z_0|$ and
$$||\,Q||(z_0,r)\ =\int\limits_{S(z_0,r)}Q(z)\ ds\ ,$$ then $f$ has a (continuous) homeomorphic
extension $\overline{f}:\overline{D}\to\overline{D'}$.}
\end{proposition}

\medskip

A domain $D\subset{\Bbb C}$ is called a {\bf quasiextremal distance
domain}, abbr. {\bf QED-domain}, see \cite{GM}, if
\begin{equation}\label{e:7.1}M(\Delta(E,F;\overline{\Bbb C})\leqslant K\cdot
M(\Delta(E,F;D))\end{equation} for some $K\geqslant1$ and all pairs
of nonintersecting continua $E$ and $F$ in $D$.

It is well known, see, e.g., Theorem 10.12 in \cite{Va}, that
\begin{equation}\label{eqKPR2.2}M(\Delta(E,F;{\Bbb C}))\geqslant\frac{2}{\pi}\log{\frac{R}{r}}\end{equation}
for any sets $E$ and $F$ in ${\Bbb C}$ intersecting all the circles
$S(z_0,\rho)$, $\rho\in(r,R)$. Hence a QED-domain has a weakly flat
boundary. One example in \cite{MRSY}, Section 3.8, shows that the
inverse conclusion is not true even among simply connected plane
domains.

A domain $D\subset{\Bbb C}$ is called a {\bf uniform domain} if each
pair of points $z_1$ and $z_2\in D$ can be joined with a rectifiable
curve $\gamma$ in $D$ such that
\begin{equation}\label{e:7.2}s(\gamma)\ \leqslant\ a\cdot|\,z_1-z_2|\end{equation}
and \begin{equation}\label{e:7.3}\min\limits_{i=1,2}\
s(\gamma(z_i,z))\ \leqslant\ b\cdot d(z,\partial D) \end{equation}
for all $z\in\gamma$ where $\gamma(z_i,z)$ is the portion of
$\gamma$ bounded by $z_i$ and $z$, see~\cite{MaSa}. It is known that
every uniform domain is a QED-domain but there exist QED-domains
that are not uniform, see \cite{GM}. Bounded convex domains and
bounded domains with smooth boundaries are simple examples of
uniform domains and, consequently, QED-domains as well as domains
with weakly flat boundaries.

\medskip

In the mapping theory and in the theory of differential equations,
it is often applied the so-called Lipschitz domains whose boundaries
are weakly flat.

\medskip

Recall first that a map $\varphi:X\to Y$ between metric spaces $X$
and $Y$ is said to be {\bf Lipschitz} provided ${\rm
dist}(\varphi(x_1),\varphi(x_2))\leqslant M\cdot{\rm dist}(x_1,x_2)$
for some $M<\infty$ and for all $x_1$ and $x_2\in X$. The map
$\varphi$ is called {\bf bi-Lipschitz} if, in addition, $M^*{\rm
dist}(x_1,x_2)\leqslant{\rm dist}(\varphi(x_1),\varphi(x_2))$ for
some $M^*>0$ and for all $x_1$ and $x_2\in X.$ Later on, $X$ and $Y$
are subsets of ${\Bbb C}$ with the Euclidean distance.

\medskip

It is said that a domain $D$ in ${\Bbb C}$ is {\bf Lipschitz} if
every point $z_0\in\partial D$ has a neighborhood $U$ that can be
mapped by a bi-Lipschitz homeomorphism $\varphi$ onto the unit disk
${\Bbb D}$ in ${\Bbb C}$ in such a way that $\varphi(\partial D\cap
U)$ is the intersection of ${\Bbb D}$ with a coordinate axes and
$f(z_0)=0$, see, e.g., \cite{Oht}. Note that bi-Lipschitz
homeomorphisms are quasiconformal and hence the Lipschitz domains
have  weakly flat boundaries.

\cc
\section{On continuous extension of lower $Q$-homeomorphisms}

\begin{lemma}\label{l:6.3} {\it
Let $D$ and $D'$ be bounded finitely connected domains in ${\Bbb C}$
and let $f:D\to D'$ be a lower $Q$-homeomorphis. If
\begin{equation}\label{e:6.4}\int\limits_{0}^{\varepsilon_0}
\frac{dr}{||Q||_{}(z_0,r)}\ =\ \infty\qquad\qquad \forall\
z_0\in\partial D\end{equation} where
$0<\varepsilon_0<d_0=\sup\limits_{z\in D}\,|z-z_0|$ and
\begin{equation}\label{e:6.6} ||Q||_{}(z_0,r)\ =\int\limits_{
S({z_0},r)}Q\ ds\ ,\end{equation} then $f$ can be extended to a
continuous mapping of $\overline{D}_P$ onto
$\overline{D^{\prime}}_P$.}
\end{lemma}

\medskip

{\bf Proof.} With no loss of generality we may assume that
$D^{\prime}$ is a circular domain and, thus,
$\overline{D^{\prime}}_P=\overline{D^{\prime}}$. By metrizability
and compactness of $\overline{D^{\prime}}$, it suffices to prove
that, for each prime end $P$ of the domain $D$, the cluster set
$$L\ =\ C(P,f)\ :=\ \left\{\ \zeta\in{{\Bbb
C}}:\ \zeta\ =\ \lim\limits_{n\to\infty}f(z_n),\ z_n\to P,\ z_n\in
D,\ n=1,2,\ldots\ \right\}$$ consists of a single point
$\zeta_0\in\partial D^{\prime}$.

Note that $L\neq\varnothing$ by compactness of the set
$\overline{D^{\prime}}$, and $L$ is a subset of $\partial
D^{\prime}$, see, e.g., Proposition 2.5 in \cite{RSal} or
Proposition 13.5 in \cite{MRSY}. Let us assume that there is at
least two points $\zeta_0$ and $\zeta_*\in L$. Set
$U=B(\zeta_0,r_0)=\{ \zeta\in\Bbb C: |\zeta -\zeta_0|<r_0\}$ where
$0<r_0<|\zeta_*-\zeta_0|$.

Let $\sigma_k$, $k=1,2,\ldots\,$, be a chain of cross--cuts of $D$
in the prime end $P$ lying on circles $S_k=S(z_0,r_k)$ from Lemma
\ref{thabc3} where $z_0 \in \partial D.$ Let $D_k$, $k=1,2,\ldots $
be the domains associated with $\sigma_k$. Then there exist points
$\zeta_k$ and $\zeta^*_k$ in the domains $D_{k}'=f(D_{k})$ such that
$|\zeta_0-\zeta_k|<r_0$ and $|\zeta_0-\zeta^*_k|>r_0$ and, moreover,
$\zeta_k\to \zeta_0$ and $\zeta^*_k\to \zeta_*$ as $k\to\infty$. Let
$C_k$ be continuous curves joining $\zeta_k$ and $\zeta^*_k$ in
$D_{k}'$. Note that by the construction $\partial U\cap
C_k\neq\varnothing$.

By the condition of strong accessibility of the point $\zeta_0$,
there is a continuum $E\subset D'$ and a number $\delta>0$ such that
$$M(\Delta(E,C_k;D'))\ \geqslant\ \delta$$ for all large enough $k$.

Without loss of generality, we may assume that the latter condition
holds for all $k=1,2,\ldots$. Note that $C=f^{-1}(E)$ is a compact
subset of $D$ and hence $\varepsilon_0={\rm dist}(z_0,C)>0$. Again,
with no loss of generality, we may assume that $r_k<\varepsilon_0$
for all $k=1,2,\ldots$.

Let $\Gamma_{m}$ be a family of all continuous curves in $D\setminus
D_m$ joining the circle $S_{0}=S(z_0,\varepsilon_0)$ and
$\overline{\sigma_m}$, $m=1,2,\ldots$. Note that by the construction
$C_k\subset D_k^{\prime}\subset D_m'$ for all $m\leqslant k$ and,
thus, by the principle of minorization
$M(f(\Gamma_{m}))\geqslant\delta$ for all $m=1,2,\ldots$.

On the other hand, the quantity $M(f(\Gamma_{m}))$ is equal to the
capacity of the condenser in $D'$ with facings $\overline{D_{m}'}$
and $\overline{f(D\setminus B_0)}$ where $B_0=B(z_0,\varepsilon_0)$,
see, e.g., \cite{Sh}. Thus, by the principle of minorization and
Theorem 3.13 in \cite{Zi}
$$M(f(\Gamma_{m}))\ \leqslant\
\frac{1}{(f(\Sigma_{m}))}$$ where $\Sigma_{m}$ is the collection of
all intersections of the domain $D$ and the spheres $S(z_0,\rho)$,
$\rho\in(r_m,\varepsilon_0)$, because
$f(\Sigma_{m})\subset\Sigma(f(S_{m}),f(S_{0}))$ where
$\Sigma(f(S_{m}),f(S_{0}))$ consists of all closed subsets of $D'$
separating $f(S_{m})$ and $f(S_{0})$. Finally, by the condition
(\ref{e:6.4}) we obtain that $M(f(\Gamma_{m}))\to0$ as $m\to\infty$.

The obtained contradiction disproves the assumption that the cluster
set $C(P,f)$ consists of more than one point. $\Box$

\cc
\section{Extension of inverse maps of lower $Q$-homeomorphisms}

The base for the proof on extending the inverse mappings of lower
$Q$-ho\-meo\-mor\-phism by prime ends in the plane is the following
fact on the cluster sets.

\bigskip

\begin{lemma}\label{l:9.1} {\it Let $D$ and $D'$ be bounded finitely
connected domains in ${\Bbb C}$, and let $f:D\to D'$ be a lower
$Q$-homeomorphism. Denote by $P_1$ and $P_2$ different prime ends of
the domain $D$ and by $\sigma_m$, $m=1,2,\ldots$, a chain of
cross--cuts in the prime end $P_1$ from Lemma \ref{thabc3}, lying on
circles $S(z_1,r_m)$, $z_1\in\partial D$, with associated domains
$d_m$. Suppose that $Q$ is integrable over $D\cap S(z_1,r)$ for a
set $E$ of numbers $r\in(0,\delta)$ of a positive linear measure
where $\delta=r_{m_0}$ and $m_0$ is a minimal number such that the
domain $d_{m_0}$ does not contain sequences of points converging to
$P_2$. If $\partial D'$ is weakly flat,
then\begin{equation}\label{e:9.2} C(P_1,f)\cap C(P_2,f)\ =\
\varnothing\ .\end{equation}} \end{lemma}

Note that in view of metrizability of the completion
$\overline{D}_P$ of the domain $D$ with prime ends, see Remark
\ref{METRIC}, the number $m_0$ in Lemma \ref{l:9.1} always exists.

\medskip

{\bf Proof.} Let us choose $\varepsilon\in(0,\delta)$ such that
$E_0:=\{r\in E:r\in(\varepsilon,\delta)\}$ has a positive linear
measure. Such a choice is possible in view of subadditivity of the
linear measure and the exhaustion $E = \cup E_m$ where $E_m = \{r\in
E:r\in(1/m,\delta)\}\,,$ $m=1,2,\ldots $. Note that by Proposition
\ref{prOS2.2}
\begin{equation}\label{e:9.3}M(f(\Sigma_{\varepsilon}))\
>\ 0\end{equation} where $\Sigma_{\varepsilon}$ is the family of all surfaces
$D(r)=D\cap S(z_1,r)$, $r\in(\varepsilon,\delta)$.

Let us assume that $C_1\cap C_2\neq\varnothing$ where
$C_i=C(P_i,f)$, $i=1,2$. By the construction there is $m_1>m_0$ such
that $\sigma_{m_1}$ lies on the circle $S(z_1,r_{m_1})$ with
$r_{m_1}<\varepsilon$. Let $d_0=d_{m_1}$ and $d_*\subseteq
D\setminus d_{m_0}$ be a domain associated with a chain of
cross--cuts in the prime end  $P_2$. Let $\zeta_0\in C_1\cap C_2$.
Choose $r_0>0$ such that $S(\zeta_0,r_0)\cap f(d_0)\neq\varnothing$
and $S(\zeta_0,r_0)\cap f(d_*)\neq\varnothing$.

Set $\Gamma=\Delta(\overline{d_0},\overline{d_*};D)$.
Correspondingly (\ref{e:9.3}), by the principle of minorization and
Theorem 3.13 in \cite{Zi},
\begin{equation}\label{e:9.4}M(f(\Gamma))\ \leqslant\
\frac{1}{M(f(\Sigma_{\varepsilon}))}\ <\ \infty\, .
\end{equation}
Let $M_0>M(f(\Gamma))$ be a finite number. By the condition of the
lemma, $\partial D'$ is weakly flat and hence there is
$r_*\in(0,r_0)$ such that
$$M(\Delta(E,F;D'))\ \geqslant\ M_0$$
for all continua $E$ and $F$ in $D'$ intersecting the circles
$S(\zeta_0,r_0)$ and $S(\zeta_0,r_*)$. However, these circles can be
joined by continuous curves $c_1$ and $c_2$ in the domains $f(d_0)$
and $f(d_*)$, correspondingly, and, in particular, for these curves
\begin{equation}\label{e:9.4a}M_0\ \leqslant\
M(\Delta(c_1,c_2;D'))\ \leqslant\ M(f(\Gamma))\,.\end{equation} The
obtained contradiction disproves the assumption that $C_1\cap
C_2\neq\varnothing$. $\Box$

\medskip

\begin{theorem}\label{thKPR8.2} {\it Let $D$ and $D'$ be bounded finitely
connected domains in ${\Bbb C}$ and $f:D\to D'$ be a lower
$Q$-homeomorphism with $Q\in L^{1}(D)$. Then $f^{-1}$ can be
extended to a continuous mapping of $\overline{D^{\prime}}_P$ onto
$\overline{D}_P$.}
\end{theorem}

\medskip

{\it Proof.} By Remark \ref{METRIC}, we may assume with no loss of
generality that $D^{\prime}$ is a circular domain,
$\overline{D^{\prime}}_P=\overline{D^{\prime}}$; $C(\zeta_0,
f^{-1})\ne\varnothing $ for every $\zeta_0\in
\partial D^{\prime}$ because $\overline{D}_P$ is metrizable and compact.
Moreover, $C(\zeta_0, f^{-1})\cap D=\varnothing $, see, e.g.,
Proposition 2.5 in \cite{RSal} or Proposition 13.5 in \cite{MRSY}.

Let us assume that there is at least two different prime ends $P_1$
and $P_2$ in $C(\zeta_0, f^{-1})$. Then $\zeta_0\in C(P_1,f)\cap
C(P_2,f)$ and, thus, (\ref{e:9.2}) does not hold. Let
$z_1\in\partial D$ be a point corresponding to $P_1$ from Lemma
\ref{thabc3}. Note that
\begin{equation}\label{eqKPR6.2ad} E\ =\ \{r\in(0,\delta):\  Q|_{D\cap S(z_1,r)}\in L^{1}(D\cap S(z_1,r))\}
\end{equation} has a positive linear measure for every $\delta>0$ by the Fubini theorem, see, e.g., \cite{Sa}, because $Q\in
L^{1}(D)$.  The obtained contradiction with Lemma \ref{l:9.1} shows
that $C(\zeta_0, f^{-1})$ contains only one prime end of $D$.

Thus, we have the extension $g$ of $f^{-1}$ to
$\overline{D^{\prime}}$ such that $C(\partial D^{\prime},
f^{-1})\subseteq \overline{D}_P\setminus D$. Really $C(\partial
D^{\prime}, f^{-1})=\overline{D}_P\setminus D$. Indeed, if $P_0$ is
a prime end of $D$, then there is a sequence $z_n$ in $D$ being
convergent to $P_0$. We may assume without loss of generality that
$z_n\to z_0\in\partial D$ and $f(z_n)\to \zeta_0\in\partial
D^{\prime}$ because $\overline{D}$ and $\overline{D^{\prime}}$ are
compact. Hence $P_0\in C(\zeta_0, f^{-1})$.

Finally, let us show that the extended mapping
$g:\overline{D^{\prime}}\to\overline{D}_P$ is continuous. Indeed,
let $\zeta_n\to\zeta_0$ in $\overline{D^{\prime}}$. If $\zeta_0\in
D^{\prime}$, then the statement is obvious. If $\zeta_0\in\partial
D^{\prime}$, then take $\zeta^*_n\in D^{\prime}$ such that
$|\zeta_n-\zeta^*_n|<1/n$ and $\rho(g(\zeta_n),g(\zeta^*_n))<1/n$
where $\rho$ is one of the metrics in Remark \ref{METRIC}. Note that
by the construction $g(\zeta^*_n)\to g(\zeta_0)$ because
$\zeta^*_n\to \zeta_0$. Consequently, $g(\zeta_n)\to g(\zeta_0)$,
too. $\Box$

\medskip

\begin{theorem}\label{t:9.12} {\it Let $D$ and $D'$ be bounded finitely
connected domains in ${\Bbb C}$. If $f:D\to D'$ is a lower
$Q$-homeomorphism with condition (\ref{e:6.4}), then $f^{-1}$ can be
extended to a continuous mapping of $\overline{D'}_P$ onto
$\overline{D}_P$.}\end{theorem}

\medskip

{\bf Proof.} Indeed, by Lemma 9.2 in \cite{KR$_1$} or Lemma 9.6 in
\cite{MRSY}, condition (\ref{e:6.4}) implies that
\begin{equation}\label{e:6.4d}\int\limits_{0}^{\delta}
\frac{dr}{||Q||_{}(z_0,r)}\ =\ \infty\qquad\qquad \forall\
z_0\in\partial D \qquad\forall\
\delta\in(0,\varepsilon_0)\end{equation} and, thus, the set
\begin{equation}\label{eqKPR6.2ad} E\ =\ \{r\in(0,\delta):\  Q|_{D\cap S(z_0,r)}\in L^{1}(D\cap S(z_0,r))\}
\end{equation}
has a positive linear measure for all $z_0\in\partial D$ and all
$\delta\in(0,\varepsilon_0)$ . The rest of arguments is perfectly
similar to one in the proof of Theorem \ref{thKPR8.2}. $\Box$
\bigskip

\cc
\section{On functions of finite mean oscillation}

Recall that a real-valued function $u$ in a domain $D$ in ${\Bbb C}$
is said to be of {\bf bounded mean oscillation} in $D$, abbr.
$u\in{\rm BMO}(D)$, or simply $u\in\ ${\bf BMO} if $u\in L_{\rm
loc}^1(D)$ and
\begin{equation}\label{lasibm_2.2}\Vert u\Vert_{*}:=
\sup\limits_{B}{\frac{1}{|B|}}\int\limits_{B}|u(z)-u_{B}|\,dm(z)<\infty\,,\end{equation}
where the supremum is taken over all discs $B$ in $D$, $dm(z)$
corresponds to the Lebesgue measure in ${\Bbb C}$ and $u_B$ is the
average of $u$ over $B$. We write $u\in{\rm BMO}_{\rm loc}(D)$ if
$u\in{\rm BMO}(U)$ for every relatively compact subdomain $U$ of $D$
(we also write BMO or ${\rm BMO}_{\rm loc }$ if it is clear from the
context what $D$ is).

\medskip

The class BMO was introduced by John and Nirenberg  in the paper
\cite{JN} and soon became an important concept in harmonic analysis,
partial differential equations and related areas, see, e.g.,
\cite{HKM} and \cite{RR}.

\medskip

A function $\varphi$ in BMO is said to have {\bf vanishing mean
oscillation} (abbreviated as $\varphi\in$ {\bf VMO}), if the
supremum in (\ref{lasibm_2.2}) taken over all disks $B$ in $D$ with
$|B| < \varepsilon$ converges to $0$ as $\varepsilon\to 0.$ VMO was
introduced by Sarason in \cite{Sarason}. Note that a large number of
papers are devoted to the existence, uniqueness and properties of
solutions for various kinds of differential equations and, in
particular, of elliptic type with coefficients of the class {VMO},
see, e.g., \cite{CFL, IwSb, MRV, Pa, Rag}.

\medskip

Following the paper \cite{IR}, see also the monographs \cite{GRSY}
and \cite{MRSY}, we say that a function $\varphi:D\to{\Bbb R}$ has
{\bf finite mean oscillation} at a point $z_0\in D$ if
\begin{equation}\label{FMO_eq2.4}\overline{\lim\limits_{\varepsilon\to0}}\ \ \
\dashint_{B(z_0,\varepsilon)}|{\varphi}(z)-\widetilde{\varphi}_{\varepsilon}(z_0)|\,dm(z)<\infty\,,\end{equation}
where \begin{equation}\label{FMO_eq2.5}
\widetilde{\varphi}_{\varepsilon}(z_0)=\dashint_{B(z_0,\varepsilon)}
{\varphi}(z)\,dm(z)\end{equation} is the mean value of the function
${\varphi}(z)$ over the disk $B(z_0,\varepsilon)$. Note that the
condition (\ref{FMO_eq2.4}) includes the assumption that $\varphi$
is integrable in some neighborhood of the point $z_0$. We say also
that a function $\varphi:D\to{\Bbb R}$ is of {\bf finite mean
oscillation in $D$}, abbr. $\varphi\in{\rm FMO}(D)$ or simply
$\varphi\in{\rm FMO}$, if $\varphi\in{\rm FMO}(z_0)$ for all points
$z_0\in D$. We write $\varphi\in{\rm FMO}(\overline{D})$ if
$\varphi$ is given in a domain $G$ in $\Bbb{C}$ such that
$\overline{D}\subset G$ and $\varphi\in{\rm FMO}(z_0)$ for all
$z_0\in\overline{D}$.

\medskip

The following statement is obvious by the triangle inequality.

\medskip

\begin{proposition}\label{FMO_pr2.1} {\it If, for a  collection of numbers
$\varphi_{\varepsilon}\in{\Bbb R}$,
$\varepsilon\in(0,\varepsilon_0]$,
\begin{equation}\label{FMO_eq2.7}\overline{\lim\limits_{\varepsilon\to0}}\ \ \
\dashint_{B(z_0,\varepsilon)}|\varphi(z)-\varphi_{\varepsilon}|\,dm(z)<\infty\,,\end{equation}
then $\varphi $ is of finite mean oscillation at $z_0$.}
\end{proposition}

In particular choosing in Proposition \ref{FMO_pr2.1},
$\varphi_{\varepsilon}\equiv0$, $\varepsilon\in(0,\varepsilon_0]$,
we obtain the following.

\medskip

\begin{corollary}\label{FMO_cor2.1} {\it If, for a point $z_0\in D$,
\begin{equation}\label{FMO_eq2.8}\overline{\lim\limits_{\varepsilon\to 0}}\ \ \
\dashint_{B(z_0,\varepsilon)}|\varphi(z)|\,dm(z)<\infty\,,
\end{equation} then $\varphi$ has finite mean oscillation at
$z_0$.} \end{corollary}

\medskip

Recall that a point $z_0\in D$ is called a {\bf Lebesgue point} of a
function $\varphi:D\to{\Bbb R}$ if $\varphi$ is integrable in a
neighborhood of $z_0$ and \begin{equation}\label{FMO_eq2.7a}
\lim\limits_{\varepsilon\to 0}\ \ \ \dashint_{B(z_0,\varepsilon)}
|\varphi(z)-\varphi(z_0)|\,dm(z)=0\,.\end{equation} It is known
that, almost every point in $D$ is a Lebesgue point for every
function $\varphi\in L^1(D)$. Thus, we have by Proposition
\ref{FMO_pr2.1} the following corollary showing that the ${\rm FMO}$
condition is natural.

\medskip

\begin{corollary}\label{FMO_cor2.7b} {\it Every
locally integrable function $\varphi:D\to{\Bbb R}$ has a finite mean
oscillation at almost every point in $D$.} \end{corollary}

\medskip

\begin{remark}\label{FMO_rmk2.13a} Note that the function $\varphi(z)=\log\left(1/|z|\right)$
belongs to BMO in the unit disk $\Bbb D$ in $\Bbb C$, see, e.g.,
\cite{RR}, p. 5, and hence also to FMO. However,
$\widetilde{\varphi}_{\varepsilon}(0)\to\infty$ as
$\varepsilon\to0$, showing that condition (\ref{FMO_eq2.8}) is only
sufficient but not necessary for a function $\varphi$ to be of
finite mean oscillation at $z_0$. Clearly, ${\rm BMO}(D)\subset{\rm
BMO}_{\rm loc}(D)\subset{\rm FMO}(D)$ and as well-known ${\rm
BMO}_{\rm loc}\subset L_{\rm loc}^p$ for all $p\in[1,\infty)$, see,
e.g., \cite{JN}. However, FMO is not a subclass of $L_{\rm loc}^p$
for any $p>1$ but only $L_{\rm loc}^1$, see examples in \cite{MRSY},
p. 211. Thus, the class FMO is essentially wider than ${\rm
BMO}_{\rm loc}$.\end{remark}

\medskip

Versions of the next lemma has been first proved for the class BMO
in the planar case in \cite{RSY$_1$} and then in the space case in
\cite{MRSY$_6$}. For the FMO case, see the papers \cite{IR} and
\cite{RS} and also the monographs \cite{GRSY} and \cite{MRSY}.

\medskip

\begin{lemma}\label{lem13.4.2} {\it Let $D$ be a domain in ${\Bbb C}$ and let
$\varphi:D\to{\Bbb R}$ be a  non-negative function  of the class
${\rm FMO}(z_0)$ for some $z_0\in D$. Then
\begin{equation}\label{eq13.4.5}\int\limits_{\varepsilon<|z-z_0|<\varepsilon_0}\frac{\varphi(z)\,dm(z)}
{\left(|z-z_0|\log\frac{1}{|z-z_0|}\right)^2}=O\left(\log\log\frac{1}{\varepsilon}\right)\
\quad\text{as}\quad\varepsilon\to 0\end{equation} for some
$\varepsilon_0\in(0,\delta_0)$ where $\delta_0=\min(e^{-e},d_0)$,
$d_0=\sup_{z\in D}|z-z_0|$.} \end{lemma}

\cc
\section{Homeomorphic extension of lower $Q$-homeomorphisms}

Combining Lemma \ref{l:6.3} and Theorem \ref{t:9.12}, we obtain the
next conclusion.

\medskip

\begin{theorem}\label{t:10.1} {\it
Let $D$ and $D'$ be bounded finitely connected domains in ${\Bbb C}$
and let $f:D\to D'$ be a lower $Q$-homeomorphism with
\begin{equation}\label{e:10.2}\int\limits_{0}^{\varepsilon_0}
\frac{dr}{||Q||_{}(z_0,r)}\ =\ \infty\qquad\qquad \forall\
z_0\in\partial D\end{equation} where
$0<\varepsilon_0<d_0=\sup\limits_{z\in D}\,|z-z_0|$ and
$$||Q||_{}(z_0,r)\ =\int\limits_{ S({z_0},r)}Q\ ds\ .$$
Then  $f$ can be extended to a homeomorphism of $\overline{D}_P$
onto $\overline{D'}_P$.}
\end{theorem}

\medskip

\begin{corollary}\label{thOSKRSS100} {\it In particular, the conclusion of Theorem
\ref{t:10.1} holds if
\begin{equation}\label{eqOSKRSS100d}q_{z_0}(r)=O\left(\log{\frac1r}\right)\qquad\qquad  \forall\
z_0\in\partial D \end{equation} as $r\to0$ where $q_{z_0}(r)$ is the
average of $Q$ over the circle $|z-z_0|=r$. }
\end{corollary}

\medskip

Using Lemma 2.2 in \cite{RS}, see also Lemma 7.4 in \cite{MRSY}, by
Theorem \ref{t:10.1} we obtain the following general lemma that, in
turn, makes possible to obtain new criteria in a great number.

\medskip

\begin{lemma}\label{lemOSKRSS1000} {\it Let $D$ and $D'$ be bounded finitely
connected domains in ${\Bbb C}$ and let $f:D\to D'$ be a lower
$Q$-homeomorphism. Suppose that
\begin{equation}\label{eqOSKRSS1000}
\int\limits_{D(z_0,\varepsilon)}Q(z)\cdot\psi_{z_0,\varepsilon}^2(|z-z_0|)\
dm(z)\ =\ o\left(I_{x_0}^2(\varepsilon)\right)\qquad\forall\
z_0\in\partial D
\end{equation} as $\varepsilon\to0$ where $D(z_0,\varepsilon)=\{z\in
D:\varepsilon<|z-z_0|<\varepsilon_0\}$ for
$0<\varepsilon_0<d(z_0)=\sup\limits_{z\in D}\,|z-z_0|$ and where
$\psi_{z_0,\varepsilon}(t): (0,\infty)\to [0,\infty]$,
$\varepsilon\in(0,\varepsilon_0)$, is a two-parameter family of
measurable functions such that
$$0\ <\ I_{z_0}(\varepsilon)\ :=\
\int\limits_{\varepsilon}^{\varepsilon_0}\psi_{z_0,\varepsilon}(t)\
dt\ <\ \infty\qquad\qquad\forall\ \varepsilon\in(0,\varepsilon_0)\
.$$ Then  $f$ can be extended to a homeomorphism of $\overline{D}_P$
onto $\overline{D'}_P$.}
\end{lemma}

\medskip

\begin{remark}\label{rmKR2.9} Note that (\ref{eqOSKRSS1000}) holds, in particular, if
\begin{equation}\label{eqOSKRSS100a} \int\limits_{B(z_0,
\varepsilon_0)}Q(z)\cdot\psi^2 (|z-z_0|)\ dm(z)\ <\
\infty\qquad\qquad \forall\ z_0\in\partial D\end{equation} where
$B(z_0,\varepsilon_0)=\{z\in D:|z-z_0|<\varepsilon_0\}$ and where
$\psi(t): (0,\infty)\to [0,\infty]$ is a measurable function such
that $I_{z_0}(\varepsilon)\to\infty$ as $\varepsilon\to0$. In other
words, for the extendability of $f$ to a homeomorphism of
$\overline{D}_P$ onto $\overline{D'}_P$, it suffices the integrals
in (\ref{eqOSKRSS100a}) to be convergent for some nonnegative
function $\psi(t)$ that is locally integrable on $(0,\infty)$ but it
has a non-integrable singularity at zero.

Note also that it is not only Lemma \ref{lemOSKRSS1000} follows from
Theorem \ref{t:10.1} but, inversely, Theorem \ref{t:10.1} follows
from Lemma \ref{lemOSKRSS1000}, too. Indeed, for the function
\begin{equation}\label{3.21}\psi_{z_0}(t)=\left \{\begin{array}{lr}
1/||Q||_{}(z_0,t),\quad & \ t\in (0,\varepsilon_0),
\\ 0,  & \ t\in [\varepsilon_0,\infty),\end{array}\right.\end{equation} we
have by the Fubini theorem that
\begin{equation}\label{3.22}\int\limits_{S(z_0,\varepsilon)} Q(z)\cdot\psi_{z_0}^2(|z-z_0|)\
dm(z)\ =\ \int\limits_{\varepsilon}^{\varepsilon_0}
\frac{dr}{||Q||_{}(z_0,r)}\ .
\end{equation}

Thus, Theorem \ref{t:10.1} is equivalent to Lemma
\ref{lemOSKRSS1000} but each of them is sometimes more convenient
for applications than another one.
\end{remark}

\bigskip

Choosing in Lemma \ref{lemOSKRSS1000} $\psi(t):=\frac{1}{t\log 1/t}$
and applying Lemma \ref{lem13.4.2}, we obtain the next result.

\medskip

\begin{theorem}\label{thOSKRSS101} {\it Let $D$ and $D'$ be bounded finitely connected domains in ${\Bbb C}$
and let $f:D\to D'$ be a lower $Q$-homeomorphism. If $Q(z)$ has
finite mean oscillation at every point $z_0\in\partial D$, then $f$
can be extended to a homeomorphism of $\overline{D}_P$ onto
$\overline{D'}_P$.}
\end{theorem}

\medskip

\begin{corollary}\label{corOSKRSS6.6.2} {\it In particular, the conslusion of Theorem \ref{thOSKRSS101} holds if
\begin{equation}\label{eqOSKRSS6.6.3}
\overline{\lim\limits_{\varepsilon\to0}}\ \
\dashint_{B(z_0,\varepsilon)}Q(z)\ dm(z)\ <\ \infty\qquad\qquad
\forall\ z_0\in\partial D\end{equation}}
\end{corollary}


\begin{corollary}\label{corOSKRSS6.6.33} {\it The conslusion of Theorem \ref{thOSKRSS101} holds if
every point $z_0\in\partial D$ is a Lebesgue point of the function
$Q:{{\Bbb C}}\to(0,\infty)$.}
\end{corollary}

\medskip

The next statement also follows from Lemma \ref{lemOSKRSS1000} under
the choice $\psi(t)=1/t.$

\medskip

\begin{theorem}\label{thOSKRSS102} {\it Let $D$ and $D'$ be
 bounded finitely connected domains in ${\Bbb C}$ and $f:D\to D'$ be
a lower $Q$-homeomorphism. If, for some
$\varepsilon_0=\varepsilon(z_0)>0$,
\begin{equation}\label{eqOSKRSS10.336a}\int\limits_{\varepsilon<|z-z_0|<\varepsilon_0}Q(z)\ \frac{dm(z)}{|z-z_0|^2}
=o\left(\left[\log\frac{1}{\varepsilon}\right]^2\right)\qquad\qquad
\forall\ z_0\in\partial D\end{equation} as $\varepsilon\to 0$, then
$f$ can be extended to a homeomorphism of $\ \overline{D}_P$ onto
$\overline{D'}_P$.}
\end{theorem}

\medskip

\begin{remark}\label{rmOSKRSS200} Choosing in Lemma \ref{lemOSKRSS1000}
the function $\psi(t)=1/(t\log 1/t)$ instead of $\psi(t)=1/t$,
(\ref{eqOSKRSS10.336a}) can be replaced by the more weak condition
\begin{equation}\label{eqOSKRSS10.336b}
\int\limits_{\varepsilon<|z-z_0|<\varepsilon_0}\frac{Q(z)\
dm(z)}{|z-z_0|\ \log{\frac{1}{|z-z_0|}}}\ =\
o\left(\left[\log\log\frac{1}{\varepsilon}\right]^2\right)\end{equation}
and (\ref{eqOSKRSS100d}) by the condition
\begin{equation}\label{eqOSKRSS10.336h} q_{z_0}(r)\ =\ o
\left(\log\frac{1}{r}\log\,\log\frac{1}{r} \right).\end{equation} Of
course, we could to give here the whole scale of the corresponding
condition of the logarithmic type using suitable functions
$\psi(t).$
\end{remark}

\medskip

Theorem \ref{t:10.1} has a magnitude of other fine consequences, for
instance:

\medskip

\begin{theorem}\label{thOSKRSS103} {\it Let $D$ and $D'$ be bounded finitely
connected domains in ${\Bbb C}$ and let $f:D\to D'$ be a lower
$Q$-homeomorphism with
\begin{equation}\label{eqOSKRSS10.36b} \int\limits_D\Phi\left(Q(z)\right)\ dm(z)\ <\ \infty\end{equation}
for a nondecreasing convex function $\Phi:[0,\infty)\to[0,\infty)$
such that
\begin{equation}\label{eqOSKRSS10.37b}
\int\limits_{\delta_*}^{\infty}\frac{d\tau}{\tau\Phi^{-1}(\tau)}=
\infty\end{equation} for $\delta_*>\Phi(0)$. Then $f$ is extended to
a homeomorphism of $\overline{D}_P$ onto $\overline{D'}_P$.}
\end{theorem}

\medskip

Indeed, by Theorem 3.1 and Corollary 3.2 in
 \cite{RSY}, (\ref{eqOSKRSS10.36b}) and
(\ref{eqOSKRSS10.37b}) imply (\ref{e:10.2}) and, thus, Theorem
\ref{thOSKRSS103} is a direct consequence of Theorem \ref{t:10.1}.

\medskip

\begin{corollary}\label{corOSKRSS6.6.3} {\it In particular, the conclusion of Theorem
\ref{thOSKRSS101} holds if
\begin{equation}\label{eqOSKRSS6.6.6}
\int\limits_{D}e^{\alpha Q(z)}\ dm(z)\ <\ \infty\end{equation} for
some $\alpha>0$.}
\end{corollary}

\medskip

\begin{remark}\label{rmOSKRSS200000}
By Theorem 2.1 in \cite{RSY}, see also Proposition 2.3 in
\cite{RS1}, (\ref{eqOSKRSS10.37b}) is equivalent to every of the
conditions from the following series:
\begin{equation}\label{eq333Y}\int\limits_{\delta}^{\infty}
H'(t)\ \frac{dt}{t}=\infty\ ,\quad\ \delta>0\ ,\end{equation}
\begin{equation}\label{eq333F}\int\limits_{\delta}^{\infty}
\frac{dH(t)}{t}=\infty\ ,\quad\ \delta>0\ ,\end{equation}
\begin{equation}\label{eq333B}
\int\limits_{\delta}^{\infty}H(t)\ \frac{dt}{t^2}=\infty\ ,\quad\
\delta>0\ ,
\end{equation}
\begin{equation}\label{eq333C}
\int\limits_{0}^{\Delta}H\left(\frac{1}{t}\right)\,dt=\infty\
,\quad\ \Delta>0\ ,
\end{equation}
\begin{equation}\label{eq333D}
\int\limits_{\delta_0}^{\infty}\frac{d\eta}{H^{-1}(\eta)}=\infty\
,\quad\ \delta_0>H(0)\ ,
\end{equation}
where
\begin{equation}\label{eq333E}
H(t)=\log\Phi(t)\ .\end{equation}

Here the integral in (\ref{eq333F}) is understood as the
Lebesgue--Stieltjes integral and the integrals in (\ref{eq333Y}) and
(\ref{eq333B})--(\ref{eq333D}) as the ordinary Lebesgue integrals.


It is necessary to give one more explanation. From the right hand
sides in the conditions (\ref{eq333Y})--(\ref{eq333D}) we have in
mind $+\infty$. If $\Phi(t)=0$ for $t\in[0,t_*]$, then
$H(t)=-\infty$ for $t\in[0,t_*]$ and we complete the definition
$H'(t)=0$ for $t\in[0,t_*]$. Note, the conditions (\ref{eq333F}) and
(\ref{eq333B}) exclude that $t_*$ belongs to the interval of
integrability because in the contrary case the left hand sides in
(\ref{eq333F}) and (\ref{eq333B}) are either equal to $-\infty$ or
indeterminate. Hence we may assume in (\ref{eq333Y})--(\ref{eq333C})
that $\delta>t_0$, correspondingly, $\Delta<1/t_0$ where
$t_0:=\sup\limits_{\Phi(t)=0}t$, set $t_0=0$ if $\Phi(0)>0$.


The most interesting of the above conditions is (\ref{eq333B}) that
can be rewritten in the form:
\begin{equation}\label{eq5!}
\int\limits_{\delta}^{\infty}\log \Phi(t)\ \ \frac{dt}{t^{2}}\ =\
\infty\ .
\end{equation}

Finally, note that if a domain $D$ in ${\Bbb C}$ is locally
connected on its boundary, then there is a natural one-to-one
correspondence between prime ends of $D$ and boundary points of $D$.
Thus, if $D$ and $D^{\prime}$ are in addition locally connected on
their boundaries in theorems of Sections 4 and 5, then $f$ is
extended to a homeomorphism of $\overline D$ onto
$\overline{D^{\prime}}$. We obtained before it similar results when
$\partial D^{\prime}$ was weakly flat which is a more strong
condition than local connectivity of $D^{\prime}$ on its boundary,
see, e.g., \cite{KPR} and \cite{KPRS}.

As known, every Jordan domain $D$ in ${\Bbb C}$ is locally connected
on its boundary, see, e.g., \cite{Wi}, p. 66. It is easy to see, the
latter implies that every bounded finitely connected domain $D$ in
${\Bbb C}$ whose boundary consists of mutually disjoint Jordan
curves and isolated points is also locally connected on its
boundary.

Inversely, every bounded finitely connected domain $D$ in ${\Bbb C}$
which is locally connected on its boundary has a boundary consisting
of mutually disjoint Jordan curves and isolated points. Indeed,
every such a domain $D$ can be mapped by a conformal mapping $f$
onto the so-called circular domain $D_*$ bounded by a finite
collection of mutually disjoint circles and isolated points, see,
e.g., Theorem V.6.2 in \cite{Goluzin}, that is extended to a
homeomorphism of $\overline D$ onto $\overline{D_*}$.

Note also that, under every homeomorphism $f$ between domains $D$
and $D^{\prime}$ in $\overline{{\Bbb C}}$, there is a natural
one-to-one correspondence between components of their boundaries
$\partial D$ and $\partial D'$, see, e.g., Lemma 5.3 in \cite{IR} or
Lemma 6.5 in \cite{MRSY}. Thus, if a bounded domain $D$ in ${\Bbb
C}$ is finitely connected and $D^{\prime}$ is bounded, then
$D^{\prime}$ is finitely connected, too.
\end{remark}

\cc
\section{Boundary behavior of mappings with finite distortion}

Recall that a homeomorphism $f$ between domains $D$ and $D'$ in
${\Bbb R}^n$, $n\geqslant2$, is called of {\bf finite distortion} if
$f\in W^{1,1}_{\rm loc}$ and \begin{equation}\label{eqOS1.3} \Vert
f'(x)\Vert^n\leqslant K(x)\cdot J_f(x)\end{equation} with some a.e.
finite function $K$ where $f^\prime(x)$ denotes the Jacobian matrix
of $f$ at $x \in D$ if it exists, $J_f(x)=\det f^\prime(x)$ is the
Jacobian of $f$ at $x$, and $\Vert f^\prime(x)\Vert$ is the operator
norm of $f^\prime(x)$, i.e.,
\begin{equation} \label{eq4.1.2}
\Vert f^\prime(x)\Vert =\max \{|f^\prime(x)h|: h \in \Rn, |h|=1\}.
\end{equation}
In the complex plane $\Vert f^\prime\Vert=|f_z|+|f_{\overline{z}}|$
and $J_f=|f_z|^2-|f_{\overline{z}}|^2$, i.e., (\ref{eqOS1.3}) is
equivalent to the condition that $K_f(z)<\infty $ a.e., see
(\ref{eq4.1.4}).

\medskip

First this notion was introduced on the plane for $f\in W^{1,2}_{\rm
loc}$ in the work \cite{IS}. Later on, this condition was replaced
by $f\in W^{1,1}_{\rm loc}$ but with the additional condition
$J_f\in L^1_{\rm loc}$ in the monograph \cite{IM}. The theory of the
mappings with finite distortion had many successors, see many
relevant references in the monographs \cite{GRSY}, \cite{HK},
\cite{IM} and \cite{MRSY}. They had as predecessors of the mappings
with bounded distortion, see \cite{Re}, and also \cite{Vo}, in other
words, the quasiregular mappings, see, e.g., \cite{HKM}, \cite{MRV},
\cite{Ri} and \cite{Vu}. They are also closely connected to the
earlier mappings with the bounded Dirichlet integral, see, e.g., the
monographs \cite{LF}, \cite{Su$_1$} and \cite{Su$_2$}, and the
mappings quasiconformal in the mean which had a rich history, see,
e.g., \cite{Ah$_4$}--\cite{Bi}, \cite{Gol$_2$}--\cite{GK},
\cite{Kr$_1$}--\cite{Ku$_1$}, \cite{Per$_1$}--\cite{Pe},
\cite{Rya$_1$}--\cite{Rya$_3$}, \cite{Str}--\cite{SS$_*$},
\cite{UV}--\cite{VU}, \cite{Zo$_1$} and \cite{Zo$_2$}.

\medskip

Note that the above additional condition $J_f\in L^1_{\rm loc}$ in
the definition of the mappings with finite distortion can be omitted
for homeomorphisms. Indeed, for each homeomorphism $f$ between
domains $D$ and $D'$ in ${\Bbb R}^n$ with the first partial
derivatives a.e. in $D$, there is a set $E$ of the Lebesgue measure
zero such that $f$ satisfies $(N)$-property by Lusin on $D\setminus
E$ and \begin{equation}\label{eqOS1.1.1}
\int\limits_{A}J_f(x)\,dm(x)=|f(A)|\end{equation} for every Borel
set $A\subset D\setminus E$, see, e.g., 3.1.4, 3.1.8 and 3.2.5 in
\cite{Fe}.

\bigskip

On the basis of Theorem \ref{thKPR3.1} and the corresponding results
on lower $Q-$ho\-meo\-mor\-phisms in Sections 5 and 6, we obtain the
following conclusions on the boundary behavior of mappings with
finite distortion.

\bigskip

\begin{theorem}\label{thKPR8.2F} {\it Let $D$ and $D'$ be bounded finitely
connected domains in ${\Bbb C}$ and $f:D\to D'$ be a homeomorphism
of finite distortion with $K_f\in L^{1}(D)$. Then $f^{-1}$ can be
extended to a continuous mapping of $\overline{D^{\prime}}_P$ onto
$\overline{D}_P$.}
\end{theorem}

\medskip

It is sufficient to assume in Theorem \ref{thKPR8.2F} that $K_f$ is
integrable only in a neighborhood of $\partial D$ and even more weak
conditions on $K_f$ due to Lemma \ref{l:9.1}.

\medskip

However, any degree of integrability of $K_f$ cannot guarantee a
continuous extension of the direct mappings $f$ to the boundary, see
an example in the proof of Proposition 6.3 in \cite{MRSY}.
Conditions for it have a perfectly different nature. The principal
related result is the following.

\medskip

\begin{theorem}\label{t:10.1F} {\it Let $D$ and $D'$ be bounded finitely
connected domains in ${\Bbb C}$, $f:D\to D'$ be a homeomorphism of
finite distortion with condition
\begin{equation}\label{ee:6.4F}\int\limits_{0}^{\varepsilon_0}
\frac{dr}{||K_f||_{}(z_0,r)}\ =\ \infty\qquad\qquad \forall\
z_0\in\partial D\end{equation} where
$0<\varepsilon_0<d_0=\sup\limits_{z\in D}\,|z-z_0|$ and
\begin{equation}\label{eq8.7.6F} ||K_f||_{}(z_0,r)=\int\limits_{
S({z_0},r)}K_f\,ds\ .\end{equation} Then  $f$ can be extended to a
homeomorphism of $\overline{D}_P$ onto $\overline{D'}_P$.}
\end{theorem}

\bigskip

Here we assume that $K_f$ is extended by zero outside of the domain
$D$.

\bigskip

\begin{corollary}\label{thOSKRSS100F} {\it In particular, the conclusion of Theorem
\ref{t:10.1F} holds if
\begin{equation}\label{eqOSKRSS100dF}k_{z_0}(r)=O\left(\log{\frac1r}\right)\ \ \ \ \ \ \forall\
z_0\in\partial D \end{equation} as $r\to0$ where $k_{z_0}(r)$ is the
average of $K_f$ over the circle $|z-z_0|=r$. }
\end{corollary}

\medskip

\begin{lemma}\label{lemOSKRSS1000F} {\it Let $D$ and $D'$ be bounded finitely
connected domains in ${\Bbb C}$ and let $f:D\to D'$ be a
homeomorphism with finite distortion. Suppose that
\begin{equation}\label{eqOSKRSS1000F}
\int\limits_{\varepsilon<|z-z_0|<\varepsilon_0}K_f(z)\cdot\psi_{z_0,\varepsilon}^2(|z-z_0|)\
dm(z)\ =\ o\left(I_{z_0}^2(\varepsilon)\right)\qquad\forall\
z_0\in\partial D
\end{equation} as $\varepsilon\to0$ where
$0<\varepsilon_0<\sup\limits_{z\in D}\,|z-z_0|$ and where
$\psi_{z_0,\varepsilon}(t): (0,\infty)\to [0,\infty]$,
$\varepsilon\in(0,\varepsilon_0)$, is a two-parameter family of
measurable functions such that
$$0\ <\ I_{z_0}(\varepsilon)\ :=\
\int\limits_{\varepsilon}^{\varepsilon_0}\psi_{z_0,\varepsilon}(t)\
dt\ <\ \infty\qquad\qquad\forall\ \varepsilon\in(0,\varepsilon_0)\
.$$ Then  $f$ can be extended to a homeomorphism of $\overline{D}_P$
onto $\overline{D'}_P$.}
\end{lemma}

\medskip

\begin{theorem}\label{thOSKRSS101F} {\it Let $D$ and $D'$ be bounded finitely
connected domains in ${\Bbb C}$ and let $f:D\to D'$ be a
homeomorphism with finite distortion. If $K_f(z)$ has finite mean
oscillation at every point $z_0\in\partial D$, then $f$ can be
extended to a homeomorphism of $\overline{D}_P$ onto
$\overline{D'}_P$.}
\end{theorem}

\medskip

In fact, here it is sufficient $K_f(z)$ to have a dominant of finite
mean oscillation in a neighborhood of every point $z_0\in\partial
D$.

\medskip

\begin{corollary}\label{corOSKRSS6.6.2F} {\it In particular, the conclusion of Theorem \ref{thOSKRSS101F} holds if
\begin{equation}\label{eqOSKRSS6.6.3F}
\overline{\lim\limits_{\varepsilon\to0}}\ \ \
\dashint_{B(z_0,\varepsilon)}K_f(z)\ dm(z)\ <\ \infty\qquad\qquad
\forall\ z_0\in\partial D\end{equation}}
\end{corollary}


\begin{theorem}\label{thOSKRSS102F} {\it Let $D$ and $D'$ be bounded finitely
connected domains in ${\Bbb C}$ and let $f:D\to D'$ be a
homeomorphism with finite distortion such that
\begin{equation}\label{eqOSKRSS10.336aF}\int\limits_{\varepsilon<|z-z_0|<\varepsilon_0}K_f(z)\
\frac{dm(z)}{|z-z_0|^2}\ =\
o\left(\left[\log\frac{1}{\varepsilon}\right]^2\right)\qquad\qquad
\forall\ z_0\in\partial D\end{equation}  Then $f$ can be extended to
a homeomorphism of $\overline{D}_P$ onto $\overline{D'}_P$.}
\end{theorem}

\medskip

\begin{remark}\label{rmOSKRSS200F} Choosing in Lemma \ref{lemOSKRSS1000}
the function $\psi(t)=1/(t\log 1/t)$ instead of $\psi(t)=1/t$,
(\ref{eqOSKRSS10.336aF}) can be replaced by the more weak condition
\begin{equation}\label{eqOSKRSS10.336bF}
\int\limits_{\varepsilon<|z-z_0|<\varepsilon_0}\frac{K_f(z)\
dm(z)}{\left(|z-z_0|\ \log{\frac{1}{|z-z_0|}}\right)^2}
=o\left(\left[\log\log\frac{1}{\varepsilon}\right]^2\right)\
.\end{equation}

Of course, we could give here the whole scale of the corresponding
conditions of the logarithmic type using suitable functions
$\psi(t).$ In particular, condition (\ref{eqOSKRSS100dF}) can be
weakened with the help of Theorem \ref{t:10.1F} by  condition
\begin{equation}\label{eqOSKRSS10.336hF} k_{z_0}(r)=O
\left(\log\frac{1}{r}\log\,\log\frac{1}{r} \right).\end{equation}
\end{remark}


Theorem \ref{t:10.1F} has a magnitude of other fine consequences,
for instance:

\medskip

\begin{theorem}\label{thOSKRSS103F} {\it Let $D$ and $D'$ be bounded finitely
connected domains in ${\Bbb C}$ and let $f:D\to D'$ be a
homeomorphism with finite distortion such that
\begin{equation}\label{eqOSKRSS10.36bF} \int\limits_D\Phi\left(K_f(z)\right)\ dm(z)\ <\ \infty\end{equation}
for a nondecreasing convex function $\Phi:[0,\infty)\to[0,\infty)$
such that
\begin{equation}\label{eqOSKRSS10.37bF}
\int\limits_{\delta_*}^{\infty}\frac{d\tau}{\tau\Phi^{-1}(\tau)}=
\infty\end{equation} for $\delta_*>\Phi(0)$. Then $f$ is extended to
a homeomorphism of $\overline{D}_P$ onto $\overline{D'}_P$.}
\end{theorem}

\medskip

\begin{corollary}\label{corOSKRSS6.6.3F} {\it In particular, the conclusion of Theorem
\ref{thOSKRSS103} holds if
\begin{equation}\label{eqOSKRSS6.6.6F}
\int\limits_{D}e^{\alpha K_f(z)}\ dm(z)\ <\ \infty\end{equation} for
some $\alpha>0$.}
\end{corollary}

\medskip

\begin{remark}\label{rmOSKRSS200000F}
Note that the condition (\ref{eqOSKRSS10.37bF}) is not only
sufficient but also necessary for a cotinuous extension to the
boundary of the mappings $f$ with integral restrictions of the form
(\ref{eqOSKRSS10.36bF}), see, e.g., Theorem 5.1 and Remark 5.1 in
\cite{KR$_3$}.
\end{remark}

\section{Boundary behavior of finitely bi--Lipschitz mappings}

Given a domain $D\subseteq\Bbb C$, following Section 5 in
\cite{KR$_7$}, see also Section 10.6 in \cite{MRSY}, we say that a
mapping $f:D\to\Bbb C$ is {\bf finitely bi-Lipschitz} if
\begin{equation}\label{eq8.12.2} 0\ <\ l(z,f)\ \leqslant\ L(z,f)\ <\
\infty \ \ \ \ \ \forall\ z\in D\end{equation} where
\begin{equation} \label{eq8.1.6} L(z,f)\ =\
\limsup_{\zeta\to z}\ \frac{|f(\zeta)-f(z)|}{|\zeta -z|}
\end{equation} and
\begin{equation}\label{eq8.1.7} l(z,f)\ =\ \liminf_{\zeta\to
z}\ \frac{|f(\zeta)-f(z)|}{|\zeta -z|}\ .\end{equation}

\medskip

The class of finitely bi-Lipschitz homeomorphisms is a natural
generalization of the well-known classes of isometries and
quasi-isometries. However, they are generally speaking are not of
finite distortion by Iwaniec.

\medskip

By the classic Stepanov theorem, see \cite{Step}, see also
\cite{Maly}, we obtain from the right hand inequality in
(\ref{eq8.12.2}) that finitely bi-Lipschitz mappings are
differentiable a.e. and from the left hand inequality in
(\ref{eq8.12.2}) that $J_f(x)\ne 0$ a.e. Moreover, such mappings
have $(N)-$property with respect to each Hausdorff measure, see,
e.g., either Lemma 5.3 in \cite{KR$_7$} or Lemma 10.6 \cite{MRSY}.
In particular, they are ACL but generally speaking are not in the
class $W^{1,1}_{\mathrm loc}$.

\medskip

However, by Corollary 5.15 in \cite{KR$_7$} and Corollary 10.10 in
\cite{MRSY}:

\medskip

\begin{lemma}\label{pr8.12.15} {\it Every finitely bi-Lipschitz
homeomorphism $f:\Omega\to\Bbb C$ is a lower $Q$-homeomorphism with
$Q=K_f$.} \end{lemma}

\medskip

\begin{corollary}\label{cor8.12.15} {\it All results on
homeomorphisms with finite distortion in Section 8 are valid for
finitely bi-Lipschitz homeomorphisms.}
\end{corollary}

\medskip

All these results for finitely bi-Lipschitz homeomorphisms are
perfectly si\-mi\-lar to the corresponding results for
homeomorphisms with finite distortion in Section 8. Hence we will
not formulate them in the explicit form here.

\bigskip

The results of this paper can be applied to the theory of boundary
value problems for the Beltrami equations including equations of the
second kind that take an important part in many problems of
mathematical physics.

\medskip
\noindent
{\bf Denis Kovtonyuk, Igor' Petkov and\\ Vladimir Ryazanov,}\\
Institute of Applied Mathematics and Mechanics,\\
National Academy of Sciences of Ukraine,\\
74 Roze Luxemburg Str., Donetsk, 83114, Ukraine,\\
denis$\underline{\ \ }$\,kovtonyuk@bk.ru, igorpetkov@list.ru,\\
vl.ryazanov1@gmail.com

\end{document}